# Primes Solutions Of Linear Diophantine Equations
## N. A. Carella

***Abstract:*** Let $k \geq 1$, $m \geq 1$ be small fixed integers, $\gcd(k, m) = 1$. This note develops some techniques for proving the existence of infinitely many primes solutions $x = p$, and $y = q$ of the linear Diophantine equation $y = mx + k$.

## 1. Introduction

Dirichlet theorem for primes in arithmetic progressions states that the linear Diophantine equation $y = mx + k \in \mathbb{Z}[x, y]$, $\gcd(k, m) = 1$, has infinitely many solutions $x = n \geq 0$ integers, and $y = p \geq 2$ primes. This note uses related analytic techniques to show that the linear Diophantine equation $y = mx + k$ has infinitely many solutions $x = p$, $y = q$, with both $p \geq 2$ and $q \geq 2$ primes.

The complementary linear Diophantine approximations problem proved in [PR] using the circle method, is quite similar. This result states that for a given tuple of arbitrary real numbers $a$, $b$, $z$, $\delta \in \mathbb{R}$ such that $a/b$ is irrational, the inequality

$$\left| ap - bq - z \right| < \delta, \tag{1}$$

has infinitely many prime solutions $p \geq 2$ and $q \geq 2$. The special case with a single prime solution $p \geq 2$ and $n \geq 1$ integer satisfies a stronger inequality $\left| ap - bn - z \right| < p^{-1/3}$, confer [HM].

The lower bounds for the numbers of prime solutions of these problems are important metrics. These two complementary problems have the same lower asymptotic formulas up to a constant. For a large number $x \geq x_0$, the asymptotic formula for number of primes solutions $p$, $q$ of the linear Diophantine equation $y = mx + k$, $\gcd(k, m) = 1$, has the lower bound $\pi_f(x) \geq bx/(\log x)^2$, $b > 0$ constant, see (18). And the asymptotic formula for the number of primes solutions $p$, $q$ of the Diophantine approximations problem has the lower bound $N(x, \delta) \geq cx/(\log x)^2$, $c > 0$ constant, see [PR].

A few other related results in the literature about primes solutions of linear equations, and configurations of primes are as follows: The result of van der Corput for three consecutive primes in arithmetic progressions states that there exist infinitely many linear polynomials $f(x) = ax + b \in \mathbb{Z}[x]$, $\gcd(a, b) = 1$, such that $f(n)$, $f(n+1)$, $f(n+2)$ are three consecutive primes for some $n \in \mathbb{Z}$. This result was recently generalized to $k$-consecutive primes $f(n)$, $f(n+1)$, ..., $f(n+k-1)$ in arithmetic progressions of arbitrary length $k \geq 1$, see [GT]. The result of Balog states that for sufficiently large $x \geq 1$, and small $k \geq 1$, there are $\geq cx^k/(\log x)^{k(k+1)/2}$ primes $k$-tuples $p_1$, $p_2$, ..., $p_k \leq x$ such that $(p_i + p_j)/2 = p$ is prime for all $i \neq j$, see [BG]. The generalization of a few of these results to integers lattices $\mathbb{Z}^d = \mathbb{Z} \times \mathbb{Z} \times \cdots \times \mathbb{Z}$ of dimensions $d \geq 1$ was announced in [CT].

The next Section has a proof of the main result, Theorem 1. This is derived from a simple primes sieve of the form $\sum_{p \leq x} \Lambda(mp + k)p^{-1}$, and a result of DeKoninck and Florian for primes in arithmetic progressions, see Theorem 3. Section 3 deals with the basic results and background materials required in the proof of Theorem 1. Section 4 presents a conditional proof of Theorem 1 based on the Montgomery conjecture for primes in arithmetic progressions, see Theorem 6, Conjecture 7 and Theorem 8.

## 2. The Main Result

The strategy of the proof of Theorem 1 is quite similar to the elementary proofs of Dirichlet theorem for primes in arithmetic progressions $\{ qn + a : n \geq 1 \}$, $\gcd(a, q) = 1$, by Selberg in [SG], Shapiro in [SP], and other authors. Other works and discussions on the elementary proofs of the prime numbers theorems are given in [DS], and [GR]. In some of these papers, it was shown that

$$\sum_{p \leq x, \, p \, \equiv \, a \bmod q} \frac{\log p}{p} > c \log x \, , \qquad (2)$$

where $c > 0$ is a constant, by elementary methods. The estimate (2) is equivalent to the fact that the logarithm derivative series

$$-\frac{1}{\varphi(q)} \sum_{\chi \bmod q} \overline{\chi}(a) \frac{L'(s, \chi)}{L(s, \chi)} = \sum_{p \, \equiv \, a \bmod q} \frac{\Lambda(n)}{n^s} \, , \qquad (3)$$

is analytic on the complex half plane $\{ s \in \mathbb{C} : \Re e(s) > 1 \}$ and has a pole at $s = 1$, see [IK, p. 37], [ST, p. 19], [TM, p. 252]. The aim of this Section is to prove that the prime pairs series $\sum_{p \geq 2} \Lambda(mp + k) p^{-s}$ for certain parameters $k, m \in \mathbb{Z}$ has similar analytic properties.

**Theorem 1.** For any fixed small parameters $k = 2a > 2$, $m = 1$, the prime pairs series $\sum_{p \geq 2} \Lambda(mp + k) p^{-s}$ is analytic on the complex half plane $\{ s \in \mathbb{C} : \Re e(s) > 1 \}$, and has a pole at $s = 1$. In particular, the Diophantine equation $y = mx + k$ has infinitely many prime pairs solutions $x = p$, $y = q$.

**Proof**: To show that the prime pairs series $\sum_{p \geq 2} \Lambda(p + 2a) p^{-s}$ has a pole at $s = 1$, consider the partial sum

$$\sum_{p \leq x} \frac{1}{p} \Lambda(p + 2a) = -\sum_{p \leq x} \frac{1}{p} \sum_{d \mid p+2a} \mu(d) \log d$$
$$= -\sum_{d \leq x+2a} \mu(d) \log d \sum_{p \leq x,\, p \equiv -2a \bmod d} \frac{1}{p}, \tag{4}$$

where $\Lambda(n) = -\sum_{d \mid n} \mu(d) \log d$, and $\gcd(2a, d) = 1$. This follows from Lemma 2 in Section 3, and inverting the order of summation.

**Remark 1:** The arithmetic progressions of small moduli $d = O(\log x)$ contribute the bulk of the main term of the final asymptotic formula (8). Moreover, the form of Theorem 3, for example, the error term has an explicit dependence on the modulo $d$, circumvents the need to use dyadic or triadic summation in the calculation of a lower bound of (4).

The totient function satisfies $d / \log d < \varphi(d) < d$, and its reciprocal inequality

$$\frac{1}{d} < \frac{1}{\varphi(d)} \leq \frac{\log d}{d}. \tag{5}$$

And the least prime $p(2a, d)$ in the arithmetic progression $\{ dn + 2a : \gcd(2a, d) = 1,$ and $n \geq 0 \}$ satisfies $d \leq dn + 2a \leq p(2a, d) < d^6$, and its reciprocal inequality

$$\frac{1}{d^6} < \frac{1}{p(2a, d)} \leq \frac{1}{dn + 2a} \leq \frac{1}{d}, \tag{6}$$

this follows from the restriction to primes in the arithmetic progression $q = p + 2a = dn + 2a$, and Linnik theorem. Moreover, there is no restriction on the range of the moduli $d \leq x$.

In view of these information, applying Theorem 3 to the inner finite sum, (with $a_0 = 2a$, $q = d$, and $\gcd(2a, d) = 1$), and simplifying yield

$$\sum_{p \leq x} \frac{1}{p} \Lambda(p + 2a) = -\sum_{\substack{d \leq x+2a, \\ \gcd(2a, d) = 1}} \mu(d) \log d \left( \frac{\log \log x}{\varphi(d)} + \frac{1}{p(2a,d)} + O\left(\frac{\log 2d}{\varphi(d)}\right) \right)$$
$$\geq -\sum_{\substack{d \leq x+2a, \\ \gcd(2a, d) = 1}} \mu(d) \log d \left( \frac{\log \log x \log d}{d} + \frac{1}{d} + O\left(\frac{\log^2 d}{d}\right) \right), \tag{7}$$

for all sufficiently large real numbers $x \geq 1$.

Apply Lemma 5, to the previous equation, for example, the first and third finite sums at $s = 1$, and $q = d$, and the second finite sum at $s = 1$, and $q = d$. This step yields

$$\sum_{p \le x+2a} \frac{1}{p} \Lambda(p+2a) \ge -\log\log x \sum_{\substack{d \le x+2a, \\ \gcd(2a,d)=1}} \frac{\mu(d)\log^2 d}{d} - \sum_{\substack{d \le x+2a, \\ \gcd(2a,d)=1}} \frac{\mu(d)\log d}{d} + O\left(\sum_{\substack{d \le x+2a, \\ \gcd(2a,d)=1}} \frac{\mu(d)\log^2 d}{d}\right)$$

$$= -\log\log x\left(-c_0 + o(1)\right) - \left(-c_1 + o(1)\right) + O\left(-c_2 + o(1)\right)$$

$$= c_0 \log\log x + o(\log\log x) ,$$

(8)

where $c_0, c_1, c_2 > 0$ are nonnegative constants, for all sufficiently large real numbers $x \ge 1$.

**Remark 2:** The error term $E(x)$ specified by the third finite sum $O(\sum \mu(d)\log^2 d / d)$ in (8) is in the range

$$c \sum_{\substack{d \le mx+k, \\ \gcd(km,d)=1}} \frac{\mu(d)\log^2 d}{d} \le E(x) \le -c \sum_{\substack{d \le mx+k, \\ \gcd(km,d)=1}} \frac{\mu(d)\log^2 d}{d} ,$$

(9)

where $c > 0$ is a constant. Thus, it is bounded by a constant, that is, $E(x) = O(1)$. This is easy to obtain using the estimate $\sum_{n \le x} \mu(n) = O(xe^{-c(\log x)^{1/2}})$, and partial summation. The trivial upper bound $E(x) = O\left(\sum_{d \le x} \log^2 d / d\right) = O(\log^3 x)$ is not used here. In fact, the trivial upper bound leads to a contradiction in (8) infinitely often as $x \to \infty$ since the finite sum

$$\sum_{\substack{d \le mx+k, \\ \gcd(km,d)=1}} \frac{\mu(d)\log^2 d}{d} = \Omega_{\pm}(x^{-1/2+\varepsilon}) ,$$

(10)

changes signs infinitely often as $x \to \infty$. More precisely, equation (7) becomes

$$2\log x \log\log x \ge \sum_{p \le x} \frac{1}{p} \Lambda(p+2a)$$

$$= -\sum_{\substack{d \le x+2a, \\ \gcd(2a,d)=1}} \mu(d)\log d\left(\frac{\log\log x}{\varphi(d)} + \frac{1}{p(2a,d)} + O(\frac{\log 2d}{\varphi(d)})\right)$$

(11)

$$\ge c_0 \log\log x + o(\log\log x) + O(\log^3 x) ,$$

which is a contradiction infinitely often as $x \to \infty$.

Now, assume that there are finitely many primes pairs $p$, and $p + 2a < x_0$, where $x_0 > 0$ is a large constant. For example, $\Lambda(p+2a) = 0$ if either $p$ or $p + 2a$ is a prime $> x_0$. Then, use a dyadic summation to compute the upper bound:

$$c_0 \log \log x + o(\log \log x) \leq \sum_{p \leq x} \frac{1}{p} \Lambda(p + 2a)$$

$$= \sum_{p \leq x_0} \frac{1}{p} \Lambda(p + 2a) + \sum_{p > x_0} \frac{1}{p} \Lambda(p + 2a) \tag{12}$$

$$\leq c_3 (\log x_0)(\log \log x_0) + c_4 + O(x^{-1/2 + \varepsilon}),$$

where $c_3, c_4 > 0$ are nonnegative constants, and $\varepsilon > 0$ is arbitrarily small, for all sufficiently large real numbers $x \geq 1$, see Lemma 4. But this is a contradiction for all sufficiently large real numbers $x > x_0$.

This shows that the partial sum $\sum_{p \leq x} \Lambda(p + 2a) p^{-1}$ of the prime pairs series $\sum_{p \geq 2} \Lambda(mp + k) p^{-s}$ is unbounded as $x \to \infty$. Therefore, it has a pole at $s = 1$, and it is analytic on the complex half plane $\{ s \in \mathbb{C} : \Re\mathrm{e}(s) > 1 \}$ since

$$\left| \sum_{p \geq 2} \frac{\Lambda(p + 2a)}{p^s} \right| \leq \sum_{n \geq 1} \frac{\Lambda(n)}{n^\sigma} < \infty, \tag{13}$$

converges for all $\Re\mathrm{e}(s) = \sigma > 1$. ∎

A closer inspection of (7) shows that it has the asymptotic formula

$$\sum_{p \leq x + 2a} \frac{1}{p} \Lambda(p + 2a) \sim c \log \log x, \tag{14}$$

where $c > 0$ is a constant. Thus, a Tauberian theorems approach can be used to derive a different proof of Theorem 1, see [KR], [RM, p. 42], and [MV, p. 261]. Moreover, observe that the better known and equivalent series

$$\sum_{n \geq 1} \frac{\Lambda(n)\Lambda(mn + k)}{n^s} \tag{15}$$

for a complex number $s \in \mathbb{C}$, seems to be more difficult to analyze, confer [KV] for related ideas.

It is an instructive exercise to show that, by partial summation, the finite sum (8) leads to the correct order of magnitude for the number of primes pairs $p, mp + k$ in the interval $[1, x]$, that is,

$$\pi_f(x) = \sum_{p,\, mp+k \,\leq\, x} 1$$

$$= \sum_{p \,\leq\, x} \frac{\Lambda(n)}{n} \frac{\Lambda(mn+k)}{\log(mn+k)} \frac{n}{\log(n)} + O(x^{1/2}) \tag{16}$$

$$\geq a \frac{x}{\log^2 x} + o(\frac{x}{\log^2 x}),$$

where $a > 0$ is a constant, and $f(x, y) = mx + k - y \in \mathbb{Z}[x, y]$. The upper bound

$$\pi_f(x) = \sum_{p,\, mp+k \,\leq\, x} 1 \leq b \frac{x}{\log^2 x}, \tag{17}$$

easily follows from (4) and Theorem 3-ii. Alternatively, it can be computed by sieve methods. In summary, this prime pairs counting function satisfies the inequality

$$a \frac{x}{\log^2 x} \leq \pi_f(x) \leq b \frac{x}{\log^2 x}. \tag{18}$$

## 3. Basic Foundation

The elementary underpinning of Theorem 1 is assembled here. The basic definitions of several number theoretical functions, and a handful of Lemmas are recorded here.

### 3.1 Formulae for the vonMangoldt Function

Let $n \in \mathbb{N} = \{\, 0, 1, 2, 3, \dots \,\}$ be an integer. The Mobius function is defined by

$$\mu(n) = \begin{cases} (-1)^v & \text{if } n = p_1 p_2 \cdots p_v, \\ 0 & \text{if } n \neq \text{squarefree}, \end{cases} \tag{19}$$

where $p_i \geq 2$ is prime. The subset of squarefree integers $\{\, n = p_1 p_2 \cdots p_v : p_i \text{ prime} \,\}$ is the support of the Mobius function $\mu : \mathbb{N} \rightarrow \{\, -1, 0, 1 \,\}$. And the vonMangoldt function is defined by

$$\Lambda(n) = \begin{cases} \log p & \text{if } n = p^k, k \geq 1, \\ 0 & \text{if } n \neq p^k, k \geq 1, \end{cases} \tag{20}$$

where $p \geq 2$ is prime. The subset of prime powers $\{\, n = p^k : p \text{ prime and } k \geq 1 \,\}$ is the support of the vonMangoldt function $\Lambda : \mathbb{N} \rightarrow \mathbb{R}$.

***Lemma* 2.**  Let $n \geq 1$ be an integer, and let $\Lambda(n)$ be the vonMangoldt function. Then

$$\Lambda(n) = -\sum_{d \mid n} \mu(d) \log d \ . \tag{21}$$

***Proof***: Use the Mobius inversion formula on the identity $\log n = \sum_{d \mid n} \Lambda(d)$ to confirm this claim.

∎

Extensive details for other identities and approximations of the vonMangoldt are discussed in [FI], [GP], [HN, p. 25], et alii.

### 3.2 Estimates For A Finite Sum Over The Primes

Effective estimates of a finite sum over the reciprocals of the primes in arithmetic progressions are provided in this Subsection.

***Theorem* 3.**  Let $x \geq 1$ be a large real number, and let $p(a, q)$ be the least prime in the arithmetic progression $\{ qn + a : n \geq 1 \}$ with $1 \leq a < q \leq x$, and $\gcd(a, q) = 1$. Then

(i)   $$\sum_{p \leq x, \, p \equiv a \bmod q} \frac{1}{p} = \frac{\log \log x}{\varphi(q)} + \frac{1}{p(a, q)} + O(\frac{\log 2q}{\varphi(q)}), \tag{22}$$

(ii)   $$\sum_{p \leq x, \, p \equiv a \bmod q} \frac{1}{p} \leq \frac{c_0}{p(a, q)} + c_1 \frac{\log \log x}{\varphi(q)},$$

where $c_0 > 0$, and $c_1 > 0$ are constants.

***Proof***: A sketch of the proof of the first part, by sieve methods, appears in [DF, p. 252, and p. 363]. The second part is derived from the Bombieri-Vinogradov inequality, and the Brun-Titchmarsh inequality.   ∎

This result simplifies some of the analysis of the distribution of primes in arithmetic progressions. Moreover, it is independent of any condition on the magnitude of the moduli; this is quite similar to the circle method analysis. In contrast, the closely related Elliot-Halberstam conjecture places a limit on the range of moduli, see [FM].

The application of Theorem 3 in Theorem 1 shows that for a properly weighted finite sum similar to (4), the sum of the error term remains smaller than the main term.

### 3.3 Estimate For A Finite Sum Over Prime Powers

Finite sums of vonMangoldt function over the primes, and prime powers occur frequently in the analysis of prime numbers. An estimate of a finite sum over the prime powers is provided in this Subsection.

***Lemma* 4.** Let $x \geq 1$ be a large real number, and let $k \geq 1$, $m \geq 1$ be small integers, $\gcd(k, m) = 1$. Then

(i) $\displaystyle\sum_{p \leq x,\, mp+k = q^v,\, v \geq 2} p^{-1} \Lambda(mp + k) \leq c_0 + O(x^{-1/2+\varepsilon}),$

(ii) $\displaystyle\sum_{p > x,\, mp+k = q^v,\, v \geq 2} p^{-1} \Lambda(mp + k) \leq c_1 + O(x^{-1/2+\varepsilon}),$ \hfill (23)

for some constants $c_0, c_1 > 0$, and $\varepsilon > 0$ is arbitrarily small.

***Proof***: Assume that $p \leq x$ is prime, and $mp + k = q^v \leq x$ is a prime power with $v \geq 2$. Let $\pi_f(x) \leq x^{1/2}$ be the prime pairs $p$ and $mp + k = q^v \leq x$ counting function. In addition, for any arbitrarily small $\varepsilon > 0$, $\Lambda(mp + k) = O(x^\varepsilon)$ as $x \to \infty$. Now, use an integral approximation to compute the estimate as follows:

$$
\begin{aligned}
\sum_{p \leq x,\, mp+k = q^v \leq x,\, v \geq 2} p^{-1} \Lambda(mp + k) &\leq c_2 \sum_{p \leq x,\, mp+k = q^v \leq x,\, v \geq 2} p^{-1+\varepsilon} \\
&\leq c_2 \int_2^x \frac{d\pi_f(t)}{t^{1-\varepsilon}} \\
&\leq c_2 \left( c_3 + O(x^{-1/2+\varepsilon}) \right),
\end{aligned}
$$
\hfill (24)

where $c_2, c_3 > 0$ are constants. ∎

### 3.4 Twisted Finite Sums

The evaluations of some common series at $s = 1$ are given in [AP], and [MV, p. 185]. Similar ideas are used here to estimate the partial sums of a few finite sums. The zeta function and the L-function are defined by $\zeta(s) = \sum_{n \geq 1} n^{-s}$, and $L(s, \chi) = \sum_{n \geq 1} \chi(n) n^{-s}$ respectively.

***Lemma* 5.** Let $s \geq 1$ be an integer, and let $x \geq 1$ be a sufficiently large real number. Then

(i) $\displaystyle\sum_{n \leq x} \frac{\mu(n) \log n}{n^s} = \frac{\zeta'(s)}{\zeta(s)^2} + O\left( x^{1-s} e^{-c(\log x)^{1/2}} \right),$ \hfill (25)

(ii) $\displaystyle\sum_{n \leq x,\, n \equiv a \bmod q} \frac{\mu(n) \log n}{n^s} = \frac{1}{q} \frac{L'(s, \chi)}{L(s, \chi)^2} + O\left( x^{1-s} e^{-c(\log x)^{1/2}} \right),$

(iii) $\displaystyle\sum_{n \leq x,\, \gcd(n,q) = 1} \frac{\mu(n) \log n}{n^s} = \frac{L'(s, \chi)}{L(s, \chi)^2} + O\left( x^{1-s} e^{-c(\log x)^{1/2}} \right),$

where $a \geq 1$, $q \geq 2$ is a pair of fixed integers, $\gcd(a, q) = 1$, $\chi$ is a character modulo $q$, and $c > 0$ is a constant.

***Proof*** (i): Write the finite sum in the form

$$\sum_{n \leq x} \frac{\mu(n) \log n}{n^s} = \sum_{n \geq 1} \frac{\mu(n) \log n}{n^s} - \sum_{n > x} \frac{\mu(n) \log n}{n^s}. \tag{26}$$

For a fixed $s \geq 1$, the constant is expressed in terms of the logarithmic derivative of the zeta function, id est,

$$\frac{d}{ds} \frac{1}{\zeta(s)} = -\sum_{n \geq 1} \frac{\mu(n) \log n}{n^s} = -\frac{\zeta'(s)}{\zeta(s)^2}. \tag{27}$$

Now, since the zeta function $\zeta(s) = \sum_{n \geq 1} n^{-s}$ is a decreasing function on the real half line $\Re e(s) > 1$, the derivative $\zeta'(s) = -\sum_{n \geq 1} (\log n) n^{-s}$ is negative on the real half line $s > 1$, see [DL, 25.2.6] for some details. The result follows from these data. The proof of (ii) is similar. To prove (iii) observe that the constant is

$$\sum_{n \geq 1, \ \gcd(n,q) = 1} \frac{\mu(n) \log n}{n^s} = -\frac{d}{ds} \frac{1}{L(s, \chi)} = \frac{L'(s, \chi)}{L(s, \chi)^2}, \tag{28}$$

where $L(s, \chi)^{-1} = \sum_{n \geq 1} \chi(n) \mu(n) n^{-s}$. ∎

A few exact evaluations of the constants are known, for example,

$$\frac{\zeta'(1)}{\zeta(1)^2} = -1 \quad \text{and} \quad \frac{L'(1, \chi)}{L(1, \chi)^2} = -\frac{\pi q}{2(h(q) \log \varepsilon)^2} \Big( \log(\sqrt{2} \Gamma(1/4)^2 / 2\pi) - (\gamma + \log 2\pi)/2 \Big), \tag{29}$$

where $h(q)$ is the class number of the quadratic field $\mathbb{Q}(\sqrt{q})$, and $\varepsilon$ is the fundamental unit. The first constant in (26) is easy to compute using the series expansion of the zeta function, see [DL, 25.2.4], and the second appears in [DM, p. 272].

## 4. Conditional Proof Of The Main Result

Let $k, m \geq 1$ be admissible integers, i. e. $\gcd(k, m) = 1$. A lower bound for the average order of the finite sum

$$\sum_{p \leq x} \Lambda(mp + k) \tag{30}$$

has a conditional derivation by means of the Montgomery conjecture for primes in arithmetic progressions, see Conjecture 7. This conjecture reduces the problem to a series of standard analytic operations involving finite sums.

***Theorem* 6.**     Assume the Montgomery conjecture. Then, for any small fixed admissible parameters $k \geq 1$, $m \geq 1$, the average order of the finite sum has the lower bound

$$\sum_{p \leq x} \Lambda(mp + k) \geq c_0 li(x) + o(li(x)),$$     (31)

where $li(x)$ is the logarithm integral, and $c_0 > 0$ is a constant. In particular, the Diophantine equation $y = mx + k$ has infinitely many prime pairs solutions $x = p$, $y = q$.

***Proof***: Suppose that $k \geq 1$ and $m \geq 1$ are small integer parameters, $\gcd(k, m) = 1$, and rewrite the finite sum over the primes as follows:

$$\begin{aligned}
\sum_{p \leq x} \Lambda(mp + k) &= -\sum_{p \leq x} \sum_{d \mid mp+k} \mu(d) \log d \\
&= -\sum_{d \leq mx+k} \mu(d) \log d \sum_{p \leq x,\, p \equiv -k/m \bmod d} 1,
\end{aligned}$$     (32)

where $\Lambda(n) = \sum_{d \mid n} \mu(d) \log d$, and $\gcd(km, d) = 1$. This follows from Lemma 2 in Section 3, and inverting the order of summation.

Applying Conjecture 7 or Theorem 8 to the inner finite sum, (with $a = -k/m$, $q = d$, and $\gcd(km, d) = 1$), and simplifying yield

$$\begin{aligned}
\sum_{p \leq x} \Lambda(mp + k) &= -\sum_{d \leq mx+k,\, \gcd(km,d)=1} \mu(d) \log d \left( \frac{li(x)}{\varphi(d)} + O(\frac{x^{1/2+\varepsilon}}{d^{1/2}}) \right) \\
&\geq -\sum_{d \leq mx+k,\, \gcd(km,d)=1} \mu(d) \log^2 d \left( \frac{li(x)}{d} + O(\frac{x^{1/2+\varepsilon}}{d^{1/2}}) \right),
\end{aligned}$$     (33)

for all sufficiently large real numbers $x \geq 1$, and an arbitrarily small number $\varepsilon > 0$.

The last inequality uses the estimate $d/\log d < \varphi(d) < d$ of the totient function. Here, there is no restriction on the range of the moduli $d \leq x$.

Apply Lemma 5-iii to the previous equation, for example, at $s = 1$, and $q = km$ in the first finite sum, and at $s = 1/2$, and $q = km$ in the second finite sum. This step yields

$$\sum_{p \leq x} \Lambda(mp + k) \geq -li(x) \sum_{d \leq mx+k,\, \gcd(km,d)=1} \frac{\mu(d) \log^2 d}{d} + O\left( x^{1/2+\varepsilon} \sum_{d \leq mx+k,\, \gcd(km,d)=} \frac{\mu(d) \log d}{d^{1/2}} \right)$$

$$\geq -li(x)\big(-c_0 + o(1)\big) + O(x^{1/2+2\varepsilon})$$

$$= c_0 li(x) + o(li(x)),$$
     (34)

where $c_0 > 0$ is a nonnegative constant, for all sufficiently large real numbers $x \geq 1$, see Remark 3 below.

Now, assume that there are finitely many primes pairs $p$, and $mp + k < x_0$, where $x_0 > 0$ is a large constant. For example, $\Lambda(mp + k) = 0$ if either $p$ or $mp + k$ is a prime $> x_0$. Then

$$
\begin{aligned}
c_0 li(x) + o(li(x)) &\leq \sum_{p \leq x} \Lambda(mp + k) \\
&= \sum_{p \leq x_0} \Lambda(mp + k) + \sum_{p > x_0} \Lambda(mp + k) \\
&\leq c_1 x_0 + \sum_{mn + k = p^v > x_0, \, v \geq 2} \Lambda(mp + k) \\
&= c_1 x_0 + O(x^{1/2 + \varepsilon}),
\end{aligned}
\tag{35}
$$

where $c_1 > 0$ is a nonnegative constant, for all sufficiently large real numbers $x \geq 1$. This follows form $\sum_{n \leq z} \Lambda(n) \leq cz$, $c > 1$ constant. But this is a contradiction for all sufficiently large real numbers $x > x_0$.

This shows that the partial sum $\sum_{p \leq x} \Lambda(mp + k)$ is unbounded as $x \to \infty$. Therefore, the prime pairs series $\sum_{p \geq 2} \Lambda(mp + k) p^{-s}$ has a pole at $s = 1$, and it is analytic on the complex half plane $\{ s \in \mathbb{C} : \Re e(s) > 1 \}$. $\blacksquare$

**Remark 3:** The error term $E(x)$ specified by the second finite sum $O(x^{1/2 + \varepsilon} \sum \mu(d) \log d / d^{1/2})$ in (31) is in the range

$$
cx^{1/2 + \varepsilon} \sum_{\substack{d \leq mx + k, \\ \gcd(km, d) = 1}} \frac{\mu(d) \log d}{d^{1/2}} \leq E(x) \leq -cx^{1/2 + \varepsilon} \sum_{\substack{d \leq mx + k, \\ \gcd(km, d) = 1}} \frac{\mu(d) \log d}{d^{1/2}},
\tag{36}
$$

where $c > 0$ is a constant. Thus, it is bounded by $E(x) = O(x^{1/2 + 2\varepsilon})$. This is easy to obtain using the estimate $\sum_{n \leq x} \mu(n) = O(x^{1/2 + \varepsilon})$, and partial summation. This estimate follows from Conjecture 7, which implies the generalized Riemann hypothesis. The trivial upper bound $E(x) = O\left(x^{1/2 + \varepsilon} \sum \log d / d^{1/2}\right) = O(x^{1 + \varepsilon})$ is not used here. In fact, it leads to a contradiction in (33) infinitely often as $x \to \infty$ since the finite sum

$$
\sum_{\substack{d \leq mx + k, \\ \gcd(km, d) = 1}} \frac{\mu(d) \log d}{d^{1/2}} = \Omega_{\pm}(x^{\varepsilon}),
\tag{37}
$$

changes signs infinitely often. More precisely, equation (33) becomes

$$x \geq \sum_{p \leq x} \Lambda(mp + k)$$

$$\geq -li(x) \sum_{d \leq mx+k,\, \gcd(km,d)=1} \frac{\mu(d) \log^2 d}{d} + O\left( x^{1/2+\varepsilon} \sum_{d \leq mx+k,\, \gcd(km,d)=} \frac{\mu(d) \log d}{d^{1/2}} \right)$$

$$\geq c_0 li(x) + o(li(x)) + O(x^{1+\varepsilon}),$$

(38)

which is a contradiction infinitely often as $x \to \infty$.

### Montgomery Conjecture

The heuristic for the Montgomery conjecture for primes in arithmetic progressions is derived from the generalized Riemann hypothesis (GRH) for $L$-functions, [MV, p. 426]. There is also a probabilistic argument similar to the Cramer model of the primes numbers. A discussion on the possible irregularity in the distribution of primes in arithmetic progressions appears in [FM].

The psi function and the prime counting functions on the arithmetic progressions $\{ qn + a : \gcd(a,q) = 1, n \geq 0 \}$ are defined by

$$\psi(x,q,a) = \sum_{n \leq x,\, n = a \bmod q} \Lambda(n) \qquad \text{and} \qquad \pi(x,q,a) = \#\{ p \leq x : p \text{ is prime} \},$$

(39)

respectively.

***Conjecture 7.*** (Montgomery 1970)   Let $\varepsilon > 0$. If $a \geq 1$, and $q \geq 2$ are integers, $\gcd(a, q) = 1$. Then

(i)  $\psi(x,q,a) = x / \varphi(q) + O(x^{1/2+\varepsilon} / q^{1/2})$,

(40)

(ii)  $\pi(x,q,a) = li(x) / \varphi(q) + O(x^{1/2+\varepsilon} / q^{1/2})$,

for all $q \leq x$, and for all large real numbers $x \geq 1$.

Below, it is shown that this conjecture is equivalent to the generalized Riemann hypothesis (GRH) for $L$- functions.

***Theorem 8.***   Let $\varepsilon > 0$. If $a \geq 1$, and $q \geq 2$ are integers, $\gcd(a, q) = 1$. Assume the GRH. Then

(i)  $\psi(x,q,a) = x / \varphi(q) + O(x^{1/2+\varepsilon} / q^{1/2})$,

(41)

(ii)  $\pi(x,q,a) = li(x) / \varphi(q) + O(x^{1/2+\varepsilon} / q^{1/2})$,

for all $q \leq x$, and for all large real numbers $x \geq 1$.

***Proof***: Suppose that the GRH holds for all the $L$-functions $L(s, \chi)$, where $\chi$ is a character modulo

$q$. By the definition of the psi function in arithmetic progressions, [IK, p. 121], [MV, p. 377], and by Theorem 9, it follows that

$$
\begin{aligned}
\psi(x,q,a) &= \frac{1}{\varphi(q)} \sum_{\chi \bmod q} \overline{\chi}(a) \psi(x,\chi) \\
&= \frac{1}{\varphi(q)} \sum_{\chi \bmod q} \overline{\chi}(a) \left( E_0(\chi)x + O(x^{1/2} \log x \log qx) \right) \\
&= \frac{x}{\varphi(q)} + O(\frac{x^{1/2} \log^2 x}{\varphi(q)}) \sum_{\chi \neq \chi_0 \bmod q} \overline{\chi}(a),
\end{aligned}
\tag{42}
$$

where $q \leq x$. Applying the Polya-Vinogradov inequality or Weil bound: $\sum_{\chi \neq \chi_0} \overline{\chi}(a) = O(q^{1/2} \log q)$, returns

$$
\psi(x,q,a) = \frac{x}{\varphi(q)} + O(\frac{x^{1/2} \log^2 x}{\varphi(q)} q^{1/2} \log q).
\tag{43}
$$

Now use $\varphi(q) \geq q/(10 \log q)$ to simplify the last expression. ∎

Note that a multiplicative character modulo $q$ has the form $\chi_u(t) = e^{i2\pi u \log t / \varphi(q)}$, where $\log t$ is the discrete logarithm modulo $q$, and $0 \leq u < \varphi(q)$. These properties make the finite sum $\sum_{\chi \neq \chi_0} \overline{\chi}(a) = \sum_u \overline{\chi}_u(a)$ of the characters a candidate for the Polya-Vinogradov inequality or Weil bound.

The twisted psi function is defined by

$$
\psi(x,\chi) = \sum_{n \leq x} \chi(n) \Lambda(n).
\tag{44}
$$

respectively.

***Theorem* 9.** ([MV, p. 425])    Let $q \geq 2$ be given, and assume the GRH for all $L$-functions $L(s, \chi)$ modulo $q$, $\gcd(a, q) = 1$. Then, for $x \geq 2$,

(i) $\psi(x,\chi) = E_0(\chi)x + O(x^{1/2} \log x \log qx)$,
$\qquad\qquad\qquad\qquad\qquad\qquad$ (45)

(ii) $\pi(x,\chi) = E_0(\chi) li(x) + O(x^{1/2} \log qx)$,

where $E_0(\chi) = 1$ if $\chi = \chi_0$, else $E_0(\chi) = 0$.

***Lemma* 10.**      Let $x \geq 1$ be a large real number and let $k \geq 1$, $m \geq 1$ be small integers, $\gcd(k, m) = 1$. Then

$$\sum_{mp+k = q^v \leq x,\, v \geq 2} \Lambda(mp + k) \leq cx^{1/2+\varepsilon},$$     (46)

for some constant $c > 0$, and $\varepsilon > 0$ arbitrarily small.

***Proof***: Assume that $p$ varies over the primes, and $q^v = mp + k$ varies over the prime powers with $v \geq 2$. Then, use the estimate $\pi(z) \leq cz / \log z$ to compute the estimate as follows:

$$\begin{aligned}
\sum_{mp+k = q^v \leq x,\, v \geq 2} \Lambda(mp + k) &\leq \log x \sum_{mp+k = q^v \leq x,\, v \geq 2} 1 \\
&\leq \log^2 x \sum_{p \leq x^{1/2}} 1 \\
&\leq \log x \left( cx^{1/2+\varepsilon} / \log x \right),
\end{aligned}$$     (47)

where $2 \leq v \leq \log x$, and $c > 1$ is a constant.                                    ∎

**Acknowledgement:** Thanks to Professor Maynard for commenting on the proof of Theorem 1.